\documentclass[a4paper,twoside,11pt]{article}

\usepackage{amsmath}
\usepackage{amsfonts}
\usepackage{amssymb}
\usepackage{amsthm}
\usepackage{amscd}

\input xy
\xyoption{all} \tolerance=500

\font\black=cmbx10 \font\sblack=cmbx7 \font\ssblack=cmbx5 \font\blackital=cmmib10  \skewchar\blackital='177
\font\sblackital=cmmib7 \skewchar\sblackital='177 \font\ssblackital=cmmib5 \skewchar\ssblackital='177
\font\sanss=cmss10 \font\ssanss=cmss8 
\font\sssanss=cmss8 scaled 600 \font\blackboard=msbm10 \font\sblackboard=msbm7 \font\ssblackboard=msbm5
\font\caligr=eusm10 \font\scaligr=eusm7 \font\sscaligr=eusm5  \font\fraktur=eufm10
\font\sfraktur=eufm7 \font\ssfraktur=eufm5

\font\bsymb=cmsy10 scaled\magstep2
\def\all#1{\setbox0=\hbox{\lower1.5pt\hbox{\bsymb
       \char"38}}\setbox1=\hbox{$_{#1}$} \box0\lower2pt\box1\;}
\def\exi#1{\setbox0=\hbox{\lower1.5pt\hbox{\bsymb \char"39}}
       \setbox1=\hbox{$_{#1}$} \box0\lower2pt\box1\;}

\newfam\bifam
\textfont\bifam=\blackital \scriptfont\bifam=\sblackital \scriptscriptfont\bifam=\ssblackital

\newfam\blfam
\textfont\blfam=\black \scriptfont\blfam=\sblack \scriptscriptfont\blfam=\ssblack

\newfam\bbfam
\textfont\bbfam=\blackboard \scriptfont\bbfam=\sblackboard \scriptscriptfont\bbfam=\ssblackboard

\newfam\ssfam
\textfont\ssfam=\sanss \scriptfont\ssfam=\ssanss \scriptscriptfont\ssfam=\sssanss
\def\sss#1{{\fam\ssfam\relax#1}}

\newfam\clfam
\textfont\clfam=\caligr \scriptfont\clfam=\scaligr \scriptscriptfont\clfam=\sscaligr

\newfam\frfam
\textfont\frfam=\fraktur \scriptfont\frfam=\sfraktur \scriptscriptfont\frfam=\ssfraktur

\def\hpb#1{\setbox0=\hbox{${#1}$}
    \copy0 \kern-\wd0 \kern.2pt \box0}
\def\vpb#1{\setbox0=\hbox{${#1}$}
    \copy0 \kern-\wd0 \raise.08pt \box0}

\def\pmb#1{\setbox0\hbox{${#1}$} \copy0 \kern-\wd0 \kern.2pt \box0}
\def\pmbb#1{\setbox0\hbox{${#1}$} \copy0 \kern-\wd0
      \kern.2pt \copy0 \kern-\wd0 \kern.2pt \box0}
\def\pmbbb#1{\setbox0\hbox{${#1}$} \copy0 \kern-\wd0
      \kern.2pt \copy0 \kern-\wd0 \kern.2pt
    \copy0 \kern-\wd0 \kern.2pt \box0}
\def\pmxb#1{\setbox0\hbox{${#1}$} \copy0 \kern-\wd0
      \kern.2pt \copy0 \kern-\wd0 \kern.2pt
      \copy0 \kern-\wd0 \kern.2pt \copy0 \kern-\wd0 \kern.2pt \box0}
\def\pmxbb#1{\setbox0\hbox{${#1}$} \copy0 \kern-\wd0 \kern.2pt
      \copy0 \kern-\wd0 \kern.2pt
      \copy0 \kern-\wd0 \kern.2pt \copy0 \kern-\wd0 \kern.2pt
      \copy0 \kern-\wd0 \kern.2pt \box0}

\def\sT{{\sss T}}

\mathchardef\za="710B  
\mathchardef\zb="710C  
\mathchardef\zg="710D  
\mathchardef\zd="710E  
\mathchardef\zve="710F 
\mathchardef\zz="7110  
\mathchardef\zh="7111  
\mathchardef\zvy="7112 
\mathchardef\zi="7113  
\mathchardef\zk="7114  
\mathchardef\zl="7115  
\mathchardef\zm="7116  
\mathchardef\zn="7117  
\mathchardef\zx="7118  
\mathchardef\zp="7119  
\mathchardef\zr="711A  
\mathchardef\zs="711B  
\mathchardef\zt="711C  
\mathchardef\zu="711D  
\mathchardef\zvf="711E 
\mathchardef\zq="711F  
\mathchardef\zc="7120  
\mathchardef\zw="7121  
\mathchardef\ze="7122  
\mathchardef\zy="7123  
\mathchardef\zf="7124  
\mathchardef\zvr="7125 
\mathchardef\zvs="7126 
\mathchardef\zf="7127  
\mathchardef\zG="7000  
\mathchardef\zD="7001  
\mathchardef\zY="7002  
\mathchardef\zL="7003  
\mathchardef\zX="7004  
\mathchardef\zP="7005  
\mathchardef\zS="7006  
\mathchardef\zU="7007  
\mathchardef\zF="7008  
\mathchardef\zW="700A  

\def\cL{{\cal L}}

\def\cM{{\cal M}}

\def\rna{\rbrack\! \rbrack}
\def\lna{\lbrack\! \lbrack}
\def\bl({\left(}
\def\br){\right)}
\def\R{\mathbb R}

\def\bt{\begin{tiny}}
\def\et{\end{tiny}}

\def\ba{\begin{array}}
\def\ea{\end{array}}

\def\ti{\times}
\def\xd{\mathrm{d}}

\def\pa{\partial}
\def\zx{\xi}
\def\ot{\otimes}

\def\Sec{\sss{Sec}}

\newcommand{\be}{\begin{equation}}
\newcommand{\ee}{\end{equation}}
\newcommand{\ra}{\rightarrow}

\newcommand{\bea}{\begin{eqnarray}}
\newcommand{\eea}{\end{eqnarray}}

\def\on{\operatorname}
\def\pa{\partial}
\def\Sec{\operatorname{Sec}}

\def\xd{\on{d}}
\def\xdp{{\on{d}^\Pi}}

\def\cM{{\mathcal M}}
\def\A{{\mathcal A}}

\def\cR{{\mathcal R}}
\def\ra{\rightarrow}
\def\we{\wedge}
\def\cL{{\mathcal L}^\Pi}

\def\ph{\phi}

\newcommand{\Ll}{\pounds}
\def\ot{\otimes}
\def\ti{\times}
\def\mod{\on{mod}}
\def\Mod{\on{Mod}}
\def\div{\on{div}}

 \def\y{\mathbf{y}}
 \def\s{\mathbf{s}}
 \def\d{\mathbf{d}}
 \def\bzx{\boldsymbol\zx}
\def\U{\mathbf{U}}
\def\uu{\mathbf{u}}

\def\cT{{\cal T}}
\def\cD{{\cal D}}
\def\half{\frac{1}{2}}

\def\*{{\textstyle *}}

\newcommand{\beas}{\begin{eqnarray*}}
\newcommand{\eeas}{\end{eqnarray*}}
\newcommand{\bfr}{\begin{frame}}
\newcommand{\efr}{\end{frame} }

\newdir{|>}{%
!/4.5pt/@{|}*:(1,-.2)@^{>}*:(1,+.2)@_{>}}
\newdir{ (}{{}*!/-5pt/@^{(}}

\newtheorem{thm}{Theorem}[section]
\newtheorem{theorem}{Theorem}[section]

\newtheorem{corollary}[thm]{Corollary}

\theoremstyle{definition}
\newtheorem{example}[thm]{Example}

\newtheorem{definition}[thm]{Definition}

\begin{document}

\markboth{Janusz Grabowski} {Modular classes revisited}

%
%

\title{MODULAR CLASSES REVISITED}

\author{JANUSZ GRABOWSKI\\ \\
Institute of Mathematics, Polish Academy of Sciences\\ \'Sniadeckich 8, 00-956 Warszawa, POLAND\\ {\tt jagrab@impan.pl}}
\date{}
\maketitle

\begin{center}
\emph{Dedicated to the memory of Jan Kubarski}
\end{center}

\medskip

\begin{abstract} We present a graded-geometric approach to modular classes of Lie algebroids and
their generalizations, introducing in this setting an idea of \emph{relative modular class} of a Dirac structure for certain type of Courant algebroids, called \emph{projectable}. This novel approach puts several concepts related to
Poisson geometry and its generalizations in a new light and simplifies proofs. It gives, in particular, a nice geometric interpretation of modular classes of twisted-Poisson structures on Lie algebroids.


\end{abstract}

\section{Introduction}

The concept of {\it modular class} of a Lie algebroid, proposed by Weinstein and introduced in \cite{ELW,We},
is a natural extension of that for a Poisson manifold \cite{Ko,We} which, in turn, can be seen as a classical reminiscence of the modular automorphism group of a von Neumann algebra,
related to quantum theories \emph{via} Tomita-Takesaki theory. It is well known that a Lie algebroid structure on a vector bundle $E$ is canonically associated with a linear Poisson structure on the dual bundle $E^\ast$, and the modular class of $E$ can be viewed as the obstruction to the existence of a homogeneous measure on $E^\ast$ being invariant with respect to all Hamiltonian vector fields on the Poisson manifold $E^\ast$ \cite{M,We}. In the framework of Lie-Rinehart algebras, the concept of modular class was developed by Huebschmann \cite{Hu}, and the role of Batalin-Vilkovisky algebras was emphasized in \cite{KS,Xu}.

Modular classes were then defined and studied also for Lie algebroid morphisms \cite{GMM,KLW,KSW}, Poisson-Nijenhuis manifolds and algebroids \cite{C,DF}, Jacobi and Poisson-Jacobi algebroids \cite{CNC}, twisted Poisson structures \cite{KSLG,SW}, twisted and quasi algebroids \cite{C1,KS1,Roy1}, even symplectic supermanifolds \cite{MV}, etc.

In \cite{Fe} Fernandes constructed a sequence of secondary characteristic classes of a Lie algebroid whose first element coincides with the modular class; see also the paper by Kubarski \cite{Ku} who, to our knowledge, was the first person defining characteristic classes and homotopies of Lie algebroids. Secondary characteristic classes of a base-preserving Lie algebroid morphism were studied also by Vaisman \cite{Vai}.

Poisson manifolds were revisited recently by Caseiro and Fernandes in \cite{CF}, where modular class of a Poisson map was defined and studied in detail. The problem was that, although any Poisson structure $\zL$ on a manifold $M$ gives rise to a Lie algebroid structure on the cotangent bundle $T^\ast M$, a Poisson map $\zf:M_1\ra M_2$ does not in general induce any canonical Lie algebroid morphism between $T^\ast M_1$ and $T^\ast M_2$. This, however, can be recognized as a particular case of a Lie algebroid relation whose modular class has been defined in \cite{G}.

On the other hand, several concepts of objects generalizing Lie algebroids appeared recently in the literature. One of them is the concept of  {\it general algebroid} or its skew-symmetric version, a {\it skew algebroid}, introduced in \cite{GU3,GU}, and applied to analytical mechanics in \cite{GG,GGU}, as an extension of the Lie algebroid/groupoid description of Lagrangian and Hamiltonian formalisms proposed by Weinstein \cite{We1} and Mart{\'\i}nez \cite{Mar}. The interest in skew algebroids, for which the Jacobi identity valid for Lie algebroids is dropped, comes also from nonholonomic mechanics in which they provide a natural geometric framework (see \cite{GG,GLMM}). The concepts of modular class and unimodularity can be extended for these, more general, algebroids almost immediately \cite{G}.

As Lie (and skew) algebroids have nice graded-geometric interpretation \cite{V}, as supermanifolds with a homological vector field, called $Q$-\emph{manifolds}. In \cite{LMS} the authors described and classified characteristic classes of $Q$-manifolds.  It is natural to ask about the description of modularity in this setting. It turns out that the modular class can be simply interpreted as the divergence of the corresponding (homological) supervector field with respect to a homogeneous Berezinian density \cite{G}.

In this note, we go further in this direction, introducing a concept of \emph{relative modular class} of a Dirac structure for a certain type of Courant algebroids which we call \emph{projectable}. This novel approach puts several concepts related to Poisson geometry and its generalizations in a new light and simplifies proofs. It gives, in particular, a nice geometric interpretation of modular classes of twisted-Poisson structures on Lie algebroids.

\section{Skew algebroids}

Let $\zt:E\to M$  be a rank-$n$ vector bundle over an $m$-dimensional manifold $M$, and  $\zp:E^*\ra M$ be its dual. We will use affine coordinates $(x^a,\zx_i)$ on $E^*$ and the dual coordinates $(x^a,y^i)$ on $E$. With $\mathcal A^i(E)=\on{Sec}(\wedge^iE)$, for $i=0,1,2,\dots$,  we denote the module of sections of the bundle $\wedge^iE$, and with $\mathcal A(E)=\bigoplus_{i\in \mathbb Z}\mathcal A^i(E)$ - the Grassmann algebra of multi-sections of $E$.

A \emph{skew algebroid} structure on $E$ can be equivalently defined as
\begin{itemize}
\item  a linear bivector
field $\zP$ on $E^*$.  In local coordinates,
\be\label{Pi} \Pi =\frac{1}{2}c^k_{ij}(x)\zx_k
\partial _{\zx_i}\wedge \partial _{\zx_j} + \zr^b_i(x) \partial _{\zx_i}
\wedge \partial _{x^b}\,,\quad\text{with}\quad c^k_{ij}(x)=-c^k_{ji}(x)\,,
\ee

\item a skew-symmetric $\R$-bilinear bracket
$[\cdot ,\cdot]_\Pi $ on the space $\Sec(E)$, together with a vector bundle morphisms\ $\zr
\colon E\rightarrow T M$ (\emph{the anchor}), such that
\be\label{qd} [X,fY]_\Pi =\zr(X)(f)Y +f [X,Y]_\Pi,
\ee
for all $f \in C^\infty (M)$, $X,Y\in \Sec(E)$,

\item a graded skew-symmetric bracket $\lna\cdot ,\cdot\rna_\Pi$
of degree $-1$, the \emph{algebroid Schouten (Gerstenhaber) bracket}, on the
Grassmann algebra $\A(E)$, satisfying the Leibniz rule
\begin{equation}\label{SchL}\lna X,Y\wedge Z\rna_\zP =\lna X,Y\rna_\zP\wedge Z+(-1)^{(k-1)l}Y\wedge\lna X,Z\rna_\zP\,,
\end{equation}
for all $X\in  \A^k(E)$,  $Y\in \A^l(E)$,

\item or as a derivation $\xdp$ of degree 1 in the Grassmann algebra
$\A(E^*)$ (\emph{de Rham derivative}). The latter is a map $\xdp:\A(E^*)\ra\A(E^*)$ such that
$\xdp:\A^i(E^*)\ra\A^{i+1}(E^*)$ and, for $\za\in\A^a(E^*)$, $\zb\in\A^b(E^*)$, we have
\be\label{dR}\xdp(\za\we\zb)=\xdp\za\we\zb+(-1)^a\za\we\xdp\zb\,.\ee
\end{itemize}
If $\Pi$ is a Poisson tensor, we speak about a \emph{Lie algebroid}.
In general, $(\xdp)^2\ne 0$, and $(\xdp)^2=0$ if and only if $\Pi$ is a Poisson tensor, thus we deal with a Lie algebroid.


For any section $X\in\Sec(E)$, the {\it Lie derivative} $\cL_X$, acting in $\A(E)$ and $\A(E^*)$,
is defined in the standard way: $\cL_X(f)=\zr(X)(f)$ for $f\in C^\infty(M)$, and
\beas
\cL_X(Y_1\we\cdots\we Y_a)&=&\sum_iY_1\we\cdots\we[X,Y_i]_\Pi\we\cdots\we Y_a\,,\\
\cL_X(\za)&=&i_X\xdp+\xdp i_X\,.\eeas

\section{Modular class of a skew algebroid}

Let now $(E,\zP)$ be a skew algebroid and $\zs=Y\ot\zm$ be a nowhere-vanishing section of the line bundle
$$L^E=\wedge^{\text{top}}E\otimes\wedge^{\text{top}}T^*M\,.$$
\begin{theorem} There is a section $\ph_\zs$ of $E^*$ (called the \emph{characteristic form of $\zs$}) such that, for all $X\in\Sec(E)$,
$$\nabla_X\zs:=\cL_X(Y)\ot\zm+Y\ot\Ll_{\zr(X)}\zm=\langle X,\ph_\zs\rangle\zs\,.$$
Moreover, if $\zs'=e^f\cdot\zs$, for some smooth function $f\in C^\infty(M)$, is another non-vanishing section, then
$\ph_{\zs'}=\ph_\zs+\xdp f$.
\end{theorem}
If the line bundle $L$ is not trivializable, we can use a nowhere-vanishing 1-density $\zs$ to define $\ph_\zs$. The class of $\ph_\zs$ in the quotient space
$[\A^1(E^*)]=\A^1(E^*)/\xdp\A^0(E^*)$ does not depend on the choice of the section $\zs$ trivializing the
bundle $L$ (resp., on the density).  We call it the \emph{modular class} of the the skew algebroid $(E,\zP)$ and denote $\mod(E,\zP)$ (orb simply $\mod(E)$ if $\Pi$ is is fixed). We call the skew algebroid \emph{unimodular} if its modular class vanishes.
Note that if $E$ is a Lie algebroid, then the characteristic form is closed,
$\xdp\ph_\zs=0$, so $\mod(E)$ is a cohomology class in $H^1(E)$.

\medskip
In local coordinates, for $\Pi$ as in (\ref{Pi}) and for the local section
$\zs=e_1\we\dots\we e_n\otimes\xd x^1\we\dots\we\xd x^m$
of $L^E$, we get
$$
\ph^\Pi_\zs=\Bigl(\sum_k c^k_{ik}(x)+\sum_a\frac{\pa\zr^a_i}{\pa x^a}(x)\Bigr)e^i.
$$
If $\Phi:E_1\to E_2$ is a morphism of skew algebroids,
$\zF^\ast\circ\xd^{\Pi_2}=\xd^{\Pi_1}\circ\,\zF^\ast$, then the \emph{modular class of $\Phi$} is defined as (cf. \cite{GMM,KSW,KLW})
$$\Mod(\Phi)=\mod(E_1)-\Phi^\*(\mod(E_2))\,.$$
 This can be generalized to skew algebroid relations $\cR\subset E_1\ti E_2$ by
\be\label{mcar}\Mod(\cR)=\pi_1^\ast(\mod(E_1))-\pi_2^\ast(\mod(E_2))\,,\ee
where $\pi_i:E_1\times E_2\supset\cR\to E_i$, $i=1,2$, is the canonical projection \cite{G}.

\section{Supergeometric description}
Let us note that, for a skew algebroid $(E,\xdp)$, we can view $\xdp$ as a (degree 1) vector field on the graded manifold $E[1]$. In local supercoordinates $(x,\y)$, associated canonically with our standard
affine coordinates $(x,y)$,
\be\label{sdR}\xdp=\frac{1}{2}c^k_{ij}(x)\y^j\y^i\pa_{\y^k} + \zr^b_i(x)\y^i\partial _{x^b}\,.
\ee
To express the modular class $\mod(E)$ directly in terms of the degree 1 vector field $\xdp$, the sheaf of sections of $L^E$ has to be replaced by the {\it
Berezinian sheaf} \ $\mathrm{Ber}=\mathrm{Ber}(E[1])$ on the supermanifold $E[1]$ whose nowhere-vanishing sections are \emph{Berezinian volumes}, representing locally \emph{Berezinian 1-desities}.

A homogeneous Berezinian volume 1-density $\s$ defines a
divergence $\div_\s$ of a homogeneous vector field $X$ by the formula (see \cite{KSM})
$$ \Ll_X\s=(-1)^{|X||s|}\s\cdot\div_\s(X)\,.$$
Here $\Ll_X\s$ is the Lie derivative of the Berezinian volume understood as a differential operator on $\A$
defined by $\Ll_X\s=-(-1)^{|X||s|}\s\circ X$.
If $\s'=e^f\s$ is another homogeneous Berezinian density on $\cD$, then
$$\div_{e^f\cdot\s}(X)=\div_\s X+X(f)\,,$$
so that the classes of $\div_\s \xdp$ and $\div_{e^f\cdot\s}(\xdp)$ coincide.
Since
$$\div_\s([X_1,X_2])=X_1(\div_\s X_2)-(-1)^{|X_1||X_2|}X_2(\div_\s X_1)\,,$$
in the case of a Lie algebroid, i.e. $[X,X]=0$ with $X=\xdp$,
we have $X(\div_\s X)=0$, so that the 1-form $\za$ represented by $\div_\s X$ is closed, so that we deal with an actual cohomology class of a Lie algebroid.
Over a superdomain $\U$ with supercoordinates $\uu=(x^1,\ldots,x^m,\y^1,\ldots,\y^n)$ the (right)
$\A(\U)$-module $\mathrm{Ber}(\U)$ is generated by
$$\s=d^{\,m|n}\uu=dx^1\wedge\ldots\wedge dx^m\otimes \pa_{\y^n}\circ\ldots\circ\pa_{\y^1}\,,$$
and the corresponding divergence of a homogeneous vector field $X=\sum_ag_a\pa_{x^a}+\sum_ih_i\pa_{\y^i}$
reads (cf. \cite{KSM}) $\div_\s(X)=\sum_a\frac{\pa g_a}{\pa x^a}-(-1)^{|X|}\sum_i\frac{\pa h_i}{\pa\y^i}$, so that
\be\label{sdiv}\div_{\s}(\xdp)=\Bigl(\sum_k c^k_{ik}(x)+\sum_a\frac{\pa\zr^a_i}{\pa x^a}(x)\Bigr)\y^i\,.
\ee
The latter superfunction represents the section $\ph^\Pi_\zs$ of $E^\ast$, so
we get the following.
\begin{theorem} Let $(E,\Pi)$ be a skew algebroid and $\xdp$ be the de Rham vector field on the supermanifold $E[1]$ corresponding to $\Pi$. Then,
the modular class $\mod(E)$ is represented by the degree-1 function $\div_{\s}(\xdp)$ being the
(super)divergence of $\xdp$ with respect to any homogeneous nowhere-vanishing Berezinian 1-density $\s$ on $E[1]$.
\end{theorem}
In what follows, we will assume that we deal with homological vector fields and Hamiltonians, i.e. Lie algebroids etc., but the major part can be directly generalized to the non-homological case of skew algebroids etc.

\subsection{Courant algebroids}
The original idea of \emph{Courant algebroid} \cite{LWX} was based on the observation
that the Vector bundle $\cT M=\sT M\oplus_M\sT^* M$, endowed with the Courant bracket\cite{Co}, plays the role of a `double'
object in the sense of Drinfeld \cite{Dr} for a pair of Lie algebroids. Let us
recall that, in complete analogy with Drinfeld's Lie bialgebras, in the category of Lie algebroids there also exist `bi-objects', Lie bialgebroids, introduced by Mackenzie and Xu \cite{MX}. On the other hand, every Lie bialgebra has a double which is a Lie algebra. This is not so for
general Lie bialgebroids. Instead, Liu, Weinstein, and Xu showed that the double of a Lie
bialgebroid is a more complicated structure they call a {\it Courant algebroid}, $\cT M$ with the Courant bracket being a special case.

In the general case, the Pontryagin bundle $\cT M$ with the
canonical symmetric pairing is replaced with a vector bundle $V\ra
M$ equipped with a nondegenerate  symmetric bilinear form $( \cdot, \cdot )$ on the bundle,
the Courant bracket is replaced with (in the Dorfman picture) a Loday (Leibniz)
bracket $[\cdot , \cdot ]$ on $\Sec (V)$, and the canonical projection $\cT M\ra\sT M$ is
replaced by a bundle map (the anchor) $\rho :V\ra TM$.

In \cite{Roy}, Roytenberg gave a nice characterization of Courant algebroids as certain Hamiltonian systems on graded symplectic manifolds.
\begin{theorem} There is a one-to-one correspondence between Courant algebroids and symplectic
N-manifolds of degree $2$, $(\cM,\zw)$, equipped with a cubic
homological Hamiltonian $H$, $\{ H,H\}=0$, where $\{\cdot,\cdot\}$ is the symplectic Poisson bracket. In this correspondence,
we identify sections of $V$ with functions of degree $1$ on
$\cM$, basic functions (functions on $M$) with functions of
degree $0$ on $\cM$, and the pseudo-Riemannian metric with the
Poisson bracket, $(X,Y)=\{ X,Y\}$. The (Dorfman) algebroid bracket
on sections of $V$ is the derived bracket $[X, Y]_H=\{\{ X,H\},Y\}$.
\end{theorem}
The Hamiltonian vector field $\d_H=\{\cdot, H\}$ is also homological, of degree 1, and defines the {\em Courant algebroid cohomology}.

If we do not assume that the cubic Hamiltonian $H$ is homological, then we deal with a \emph{nonhomological Courant algebroid} whose derived bracket on sections of $V$ is still well defined but does not satisfies the Jacobi identity.

Consider local coordinates $(x^a,\zz^i,p_b)$ in $\cM$,
corresponding to coordinates $(x^a)$ on $M$ and a local basis $\{
e_i\}$ of sections of $V$ such that $(e_i,e_j)=g_{ij}=const$\,,
$e_i=g_{ij}\zz^j$ interpreted as linear functions on $V$. Then, the symplectic
form $\zw$ reads \be\zw=\xd p_a\xd x^a+\frac{1}{2}g_{ij}\xd\zz^i
\xd\zz^j\,,\ee and any cubic Hamiltonian is of the form \be
H=\zz^i\zr^a_i(x)p_a-\frac{1}{6}\zvf_{ijk}(x)\zz^i\zz^j\zz^k\,.\ee
For the corresponding Courant algebroid, the Dorfman bracket and
the anchor are uniquely determined by
$$([e_i,
e_j]_H,e_k)=\zvf_{ijk}(x)\,,\quad
\rho(e_i)=\zr^a_i(x)\partial_{x^a}\,.
$$

\section{Relative modular class of a Dirac structure}
We can try to define the modular class of a Courant algebroid $(\cM,\zw,H)$ similarly as above, as the divergence of the Hamiltonian vector field $\d_H$ with respect to a homogeneous Berezinian density $\s$ on $\cM$. Changing the density by a positive factor $e^f$ of degree 0 results in changing the divergence by $\d_H(f)$ which is cohomologically irrelevant. Unfortunately, it is easy to see that we can always find a homogeneous `Liouville density' $\s$ for which $\div_{\s}(\d_H)=0$, so that this class is always trivial.

However, the divergence of a Hamiltonian vector field need not be trivial on a submanifold. This is e.g. the case of a Lie algebroid $E$, considered as a graded submanifold $E[1]$ of of the Courant algebroid $\cM_E=\sT^\*[2]E[1]\simeq\sT^\*[2]E^\*[1]$ with the Hamiltonian expressed in the natural coordinates $(x,\y,\bzx,p)$ as
$$H=\frac{1}{2}c^k_{ij}(x)\y^j\y^i\bzx_k + \zr^b_i(x)\y^ip_b\,.$$
Here, $V=E\oplus_ME^\*$ has a canonical symmetric pairing.
The corresponding Hamiltonian vector field, restricted to $E[1]$, is (\ref{sdR}) and its divergence represents $\mod(E)$.

The above example can be immediately extended to {\em Dirac structures}.
In the graded-geometrical language (see e.g. \cite{Gru,Sev}), they are simply lagrangian submanifolds $\cD$ of a Courant algebroid $(\sT^\*[2]E[1],\zw_{E[1]},H)$ such that $\d_H$ is tangent to $\cD$. In what follows, Courant algebroids will be always of the form $(\sT^\*[2]E[1],\zw_{E[1]},H)$, denoted with $(\cM_E,H)$, therefore they are canonically bi-graded. Dirac structures correspond to maximally isotropic vector subbundles $\cD_E$ of $E\oplus_ME^\*$ (supported on a submanifold $\cD_M$ of $M$), closed with respect to the derived Dorfman bracket. This bracket defines actually a Lie algebroid structure on $\cD_E$. The divergence $\za=\div_\s \d_\cD$ of the restriction to $\cD$ of the Hamiltonian vector field $\d_H$,  $\d_\cD=(\d_H)_{|\cD}$, relative to a homogeneous Berezinian density $\s$ on $\cD$, represents the modular class $\mod(\cD)$ of the Lie algebroid $(\cD, \d_\cD)$.

\medskip
Let us assume now that the Hamiltonian vector field $\d_H$ is projectable with respect to the canonical projection $\pi_E$ of $\sT^\*[2]E[1]$ onto $E[1]$. It is easy to see that this corresponds to the fact that  in the decomposition $H=\zm+\zg+\zvf+\psi$ of the cubic Hamiltonian into the parts of bi-degrees $(2,1), (1,2), (0,3)$, and $(3,0)$, respectively, the parts $\zg$ and $\psi$ vanish. This is a particular case of a \emph{Lie quasi-bialgebroid} in the terminology of \cite{KS1}. The corresponding Courant algebroid we will call \emph{projectable}.

\begin{theorem} In any projectable Courant algebroid $(\cM_E,H=\zm+\zvf)$, the projection $\d_E=(\pi_E)_*\d_H=(\pi_E)_*\d_\zm$ of the Hamiltonian vector field $\d_H$ is a homological vector field of degree 1 on $E[1]$, therefore defining a Lie algebroid structure.
Moreover, $\d_E$ and $\d_\cD$ are $(\pi_E)_{|\cD}$ related for any Dirac structure $\cD\subset\cM_E$, so that
\be\label{rmcD}\pi_E^\cD=(\pi_E)_{|\cD_E}:\cD_E\to E
\ee
is a Lie algebroid morphism.
\end{theorem}
\begin{definition}
The {\em relative modular class} $\Mod(\cD)$ of a Dirac structure $\cD$ in a projectable Courant algebroid $(\cM_E,H=\zm+\zvf)$ is the modular class $\Mod(\pi_E^\cD)$ of the Lie algebroid morphism (\ref{rmcD}).
\end{definition}
\begin{example} Let us see that the modular class of a \emph{twisted Poisson structure} \cite{Roy1,SW}, as defined in \cite{KSLG}, is a particular case of the above.
For a projectable Courant algebroid $(\cM_E,H=\zm+\zvf)$, consider a `bivector field' $P\in\Sec(\we^2 E)$ and its graph $\cD^P\subset E\oplus_M E^\*$,
$\cD^P_E=\{ (P^\natural(\za),\za):\za\in E^\*\}$, which is clearly a maximally isotropic subbundle in $E\oplus_M E^\*$, so corresponds to a lagrangian submanifold $\cD^P\subset\cM_E$. It defines a Dirac structure if and only if $\d_H$ is tangent to $\cD^P$, i.e. the cubic Hamiltonian $H=\zm+\zvf$ is constant, thus 0, on $\cD^P$. Interpreting $P$ as a quadratic function on $E^\*$, $P=\half P^{ij}\bzx_i\bzx_j$, we get $\cD^P$ as the image of the lagrangian submanifold $E^\*[1]\subset\cM_E$ under the time 1 map  $\exp X_P$ of the flow of the Hamiltonian vector field $X_P=\{\cdot,P\}$. A function $F$ vanishes on this image if and only if its pull-back
$$(\exp(-X_P))^\*F=F-X_P(F)+\half X^2_P(F)+\cdots=F+\{P,F\}+\half\{P,\{P,F\}\}+\cdots$$
vanishes on $E^\*[1]$. Hence, $(\exp(-X_P))^\*H$ vanishes on $E^\*[1]$ if and only if
$$\half\{P,\{P,\zm\}\}+\frac{1}{6}\{P,\{P,\{P,\zvf\}\}\}=-\half\lna P,P\rna_{\d_E}+(\we^3P^\natural)\zvf=0\,,$$
where $\lna\cdot,\cdot\rna_{\d_E}$ is the Schouten bracket of the Lie algebroid $(E,\d_E)$ and $\zvf$ is interpreted as a section of $\we^3E^\*$. This means exactly that $(P,\zvf)$ is a quasi-Poisson structure of the Lie algebroid $(E,\d_E)$.
Identifying $\cD^P$ with $E^\*$ \emph{via} $\exp X_P$, the Lie algebroid structure on $E^\*$ coming from $(\cD,\d_\cD)$ is the one defined by the projection of the Hamiltonian vector field $\d_{E^\*}^{P,\zvf}$ with the Hamiltonian $\{P,\zm\}+\half\{P,\{P,\zvf\}\}$. The corresponding Lie bracket on sections of $E^\*$ reads
$$[\za,\zb]_{P,\zvf}=\Ll_{P^\natural\za}\zb-\Ll_{P^\natural\zb}\za-\xd(P(\za,\zb))+\zvf(P^\natural\za,P^\natural\zb,\cdot)\,.$$ This is the Lie algebroid associated with a quasi-Poisson structure described in \cite{SW} and \cite{Roy1}. Since, under the identification $\cD^P_E\simeq E^*$, the projection $\pi_E^\cD$ corresponds to $P^\natural$,
we get for free the following result whose known up to date proofs were rather technical.
\begin{corollary}
The map $P^\natural:(E^*,[\cdot,\cdot]_{P,\zvf})\to (E,[\cdot,\cdot]_{\zm})$ is a Lie algebroid morphism.
\end{corollary}
\noindent Moreover,
the relative modular class $\Mod(\cD^P)$ corresponds to
$$\mod(E^\*)-(P^\natural)^\*(\mod(E))=\mod(E^\*)+P^\natural(\mod(E))\,,$$
which is the modular class of the quasi-Poisson structure $(P,\zvf)$ associated with the Lie algebroid $(E,\d_E)$ as it appears in \cite{KSW}.
\end{example}
Note finally that our approach \emph{via} Dirac structures of projectable Courant algebroids can also be applied to other `twisted' structures, like \emph{quasi-Poisson $G$-manifolds} in the sense of \cite{AKS,AKSM,BCS}, or `twisted Nambu structures' in the context of \emph{higher Courant algebroids} \cite{BS,Za}.

Another generalization depends on considering nonhomological Courant algebroids whith relations to \emph{quasi/almost-Poisson brackets} used in the nonholonomic mechanics \cite{GG,GLMM,GU,vSM}.
\section{Acknowledgment}
Research  founded by the  Polish National Science Centre grant under the contract number DEC-2012/06/A/ST1/00256.

\end{document}